\def\R{{\bf R}}
\magnification = 1200

\centerline {\bf Rectangle Groups}
\medskip
\centerline {M. J. Dunwoody}

\bigskip
\noindent
{\bf Introduction}

We describe a class of groups which illustrate various themes in geometric 
group theory.
Each of these groups contains a maximal free abelian subgroup.  In fact each one  contains 
infinitely many such subgroups.
The actions of free abelian groups on $\R$-trees is well understood.
If $F$ is a free abelian group with generating set $a_1, a_2, \dots , a_n$ then
for any choice of real numbers $r_1, r_2, \dots , r_n$ there is a homomorphism
$\theta : F \rightarrow \R$ in which $\theta (a_i) = r_i$.   Here both groups
are written additively.   We can regard $\R$ as an $\R$-tree on which $\R$ acts via translation 
and $F$ acts
via the homomorphism $\theta $.  Any action of $F$ on an $\R$-tree is, essentially, like such
an action.
We will see that for each of our new groups $G$ an action of $G$ on an $\R$-tree
corresponds uniquely to an  action of the maximal free abelian subgroup $F$ and there is an action
of $G$ for each such action of $F$.

In the theory of JSJ-decompositions of a finitely presented group $G$  for splittings over slender
subgroups the group $G$ has a decomposition as a fundamental group of a graph of groups
and incompatible decompostions of $G$ over slender subgroups correspond to geometric
splittings of vertex groups of this decomposition.   Some  vertex groups   map to  $2$-orbifold groups and there are incompatible decompositions corresponding to curves which intersect  essentially in the orbifold.  Any other vertex group has 
a homomorphic image which is a group of isometries of $n$-dimensional Euclidean space for $n\geq 3$.
If one was attempting to generalize JSJ-theory to splittings over non-small subgroups then
the groups above would have to occur as vertex groups in the decompositions.  They do not map to orbifold groups and they
do not map onto groups of isometries of Euclidean space, though they do have subgroups with this 
property.

I thank Cormac Long for helping me by investigating these groups with GAP.

I am very grateful to Gilbert Levitt for correcting some of the wilder ideas I had about these groups.
\medskip
\noindent
{\bf The groups}

Let $R$ be the group with  presentation

$$R= \langle a, b, c, d  \  |  \  a^{n_a} = b^{n_b} = c^{n_c} = d^{n_d} = 1,
ba^{-1} = dc^{-1}, ca^{-1} = db^{-1}\rangle .$$

Here $n_a, n_b, n_c, n_d$ are integers $\geq 2$.

The form of this presentation that suggests calling $R$ a {\it rectangle } group.
I originally called these groups parallelogram groups but after further study if seems
that in most  geometric interpretations of these groups   the angles involved are right angles.

First note that $R$ has two incompatible decompositions as a free product
with amalgamation,

$$R =   \langle a, b \rangle *_C \langle c, d \rangle $$

where $C = \langle  a{b}^{-1}, {b}^{-1}a \rangle = \langle cd^{-1}, d^{-1}c \rangle $.

or

$$R =   \langle a, c\rangle *_C \langle b, d\rangle $$

where $C = \langle  ac^{-1}, c^{-1}a \rangle = \langle bd^{-1}, d^{-1}b \rangle $.

In these decompositions $C$ is free of rank two unless at least two of $n_a, n_b, n_c, n_d$
are two.

Note that $x = ab^{-1} = cd^{-1}$ commutes with $y = ac^{-1} = bd^{-1}$ since
$xy = yx = ad^{-1}$, and $x, y$ generate a free abelian group of rank $2$.
This subgroup is denoted $F$.

    \bigskip
    
    If $n_a = n_b = n_c = n_d = 2$, then $C = <ab> = <cd>$ is infinite cyclic and $R$ is the free product
    with amalgamation of two infinite dihedral groups.
    In this case $R$ has an action on the Euclidean plane.
    Let $A, B, C, D$ be the vertices of a rectangle as in Fig 1.
    Let $a, b, c, d$ be rotations through $\pi$ at $A, B, C, D$ respectively.
    Then this gives the action of $R$ on the plane.   Note that $x = ab^{-1}$ will be a translation of two units
    in the horizontal  direction and $y = ac^{-1}$ will be translation by two units in a vertical direction.
    In the action of $R$ there will be two orbits of rectangles giving a chequer board pattern as shown.
    There will be $4$ orbits of both vertices and edges.  One gets a transversal for the action
    on the $1$-skeleton by taking the edges and vertices of a single rectangle.
    
    For the general situation in which $n_a, n_b, n_c, n_d$ are not all two, we construct
    a simply connected $2$-dimensional $R$-complex $X$ in which again there are $4$ orbits of vertices and edges
    and there is a sub-complex which is a plane on which $x, y$ act as above.
    We will in fact only consider the case when each of $n_a, n_b, n_c, n_d$ is at least $3$, but there is
    a similar theory for the omitted cases.
    Let $A, B, C, D$ be a tranversal for the vertex set of $X$, and let the stabilizers of $A, B, C, D$ be generated by $a, b, c, d$ respectively.
    Let $e = \vec AB, f = \vec AC, g = \vec BD, h = \vec CD$ be a tranversal for the action on the directed edges of $X$ each of
    which has trivial stabilizer.  Clearly $A, B, D, C$ are the vertices of a cycle with edges
    $e, g, \bar h, \bar f$
    However so also are $A, aB, cD, C$ since $aB =ab^{-1}B$ and $cD =cd^{-1}D$ and so 
    the edge $g = BD$ acted on by $ab^{-1} = cd^{-1}$ is the edge with vertices $aB, cD$.
    Similarly there is an edge joining $aC$ and $bD$ and so there is a cycle $A, aC, bD, B$.
    Note that we are assuming the $1$-skeleton of $X$ is a simple graph, so that there is only
    one directed edge with a given ordered pair of vertices.
    The three cycles $(A, B, D, C)$, $(A, aB, cD, C)$, $(A, aC, bD, B)$ are in different $R$-orbits.
    In Fig 1 there is another cycle containing $A$, namely $A, aB, aD, aC$.  Clearly
    this is in the same orbit as $(A, B, D, C)$
    We can therefore attach three $R$-orbits of $2$-cells and obtain the required complex $X$, in which
    the $(x, y)$ plane is indeed a plane.  To summarize,  in Fig 1 the shaded rectangles containing 
    $A$ as a boundary point are in the same orbit, but the unshaded rectangles containing $A$ as
    a boundary point are in different orbits.
    A vertical line in the $(x, y)$-plane $E$ will be part of a track in $X$ which will correspond to one of
    the splittings of $R$ over a free rank $2$ subgroup.   This track will be a $3$-regular tree
    and its stabilizer in $R$ will have $2$ orbits of vertices.  Here we note that $y$ is a translation
    by $2$ vertical units.  There is no element in $R$ translating through one vertical unit.
    A horizontal line in the $(x, y)$-plane is part of a track giving the other splitting of $R$.
    The space $X$ is simply connected.  It is the universal cover of a developable complex of groups
    which is shown to be simply connected in [H].
    One can also establish that $X$ is simply connected by cutting it up along the simply connected 
    tracks just discussed and considering the pieces that remain. 
 
    The link of the vertex $A$ in $X$ will be a $3$-regular bipartite graph with $2n_a$ vertices.
    The vertices of the link correspond to  the edges $(A, gB)\ $ for $  g \in \langle a \rangle $ and 
    $(A, g'C)  $ for $ g'  \in \langle a \rangle $ and the edges  to the $2$-cells which include the three
    vertices $g'C, A, gB$.
    
      Since $X$ is locally $CAT(0)$ - since the  link
   of any vertex has no cycles of length less than $4$ - and  it is simply connected, it will be
    $CAT(0)$ by [BH, p206].
    
    Another - probably better - way of constructing $X$ is to use the cubing construction of [Sa].
    This is the approach that can best be generalized in higher dimensional generalization.
    The two decompositions for $R$ we described earlier give two multi-ended pairs $(R, C_1), (R, C_2)$
    where $C_1  = \langle  ab^{-1}, b^{-1}a \rangle = \langle cd^{-1}, d^{-1}c \rangle $
    and $C_2 = \langle  ac^{-1}, c^{-1}a \rangle = \langle bd^{-1}, d^{-1}b \rangle $.
    Using such data Sageev describes how to construct a $1$-connected cube complex which will be our
    space $X$.  The two decompositions are associated with hyperplanes in this complex.

    Since $X$ is simply connected `tracks' separate (see [DD]).  Here we take a track to
    be a connected subspace that intersects  any $2$-cell in finitely many parallel straight lines that do not           
    intersect
    a vertex.  
    
    A  surjective homomorphism from $F$ to ${\bf Z}$  the integers under addition is determined
    by a pair $(p, q)$ of coprime integers.  One way of getting tracks in the $(x, y)$-plane corresponding
    to such a homomorphism, in the case when $p, q \not= 0$,  is to assign length $p$ to the edge $AB$ and length $q$ to the edge $AC$
    and then foliate the plane by lines of gradient one (See Fig 2), and which intersect the axes at points with integer coefficients.  Each  leaf in $E$ maps to one of two  
    $(p,q)$-curves in the torus $F\backslash E$.   Note that the action of $x$ is translation through
    $2p$ and the action of $y$ is translation through $2q$.   
    In the quotient space $R\backslash X$ the leaves maps to a track which has $2(p+q)$ intersections
    with the $1$-skeleton.   It has $p$ intersection points with each of the edges $AB$ and $CD$ and
    $q$ intersections with each of $AC$ and $BD$.   In $X$ this track will lift to an $R$-pattern of tracks.
    Each component track will be a $3$-regular tree.  The stabilizer of this track will be a free group.
    The quotient graph has $2(p+q)$ vertices and - since it is a $3$-regular graph - it has $3(p+q)$ edges.
    The rank of the stabilizer will therefore be $p+q+1$.
    Apart from the splittings described above, which can be regarded as corresponding to the pairs 
    $(1,0), (0,1)$, there is a splitting of $R$ for each pair of coprime positive integers $(p, q)$.

\bigskip
    \noindent
    {\bf The tree}

    If we assign positive lengths  $\alpha , \beta $ to the edges of $X$ so that
    $|AB| = |CD| = \alpha , |AC| = |BD| = \beta $, then we obtain a metric on $X$   in which each $2$-cell    has  the topology of a rectangle in the Euclidean plane with sides
    $\alpha, \beta$. 
     There will be a foliation on $X$ invariant under $R$ in which the $(x, y)$ plane is foliated
    by straight lines with gradient one as in Fig 2.   If $\alpha, \beta $ are dependent over the rationals
    then the image of a leaf in $X$ which does not intersect a vertex will be a track in $R\backslash X$
    and will correspond to one of the decompositions of $R$ given earlier, where $\alpha = p\delta,
    \beta = q\delta $ where $p, q $ are coprime integers.
    We now investigate what happens if  $\alpha , \beta $ are independent over the rationals
    The action of $R$ on a simplicial tree as described above will determine an element of 
    PLF$(R)$ the projectivized length functions on $R$, which is a compact space (see [Sh]).
    By approximating $\alpha , \beta $ by rationals we obtain a sequence of actions on simplicial
    trees which will tend to a non-simplicial action of $R$ on an $\R$-tree $T$,  which restricts to the 
    obvious action of $F$ on $\R$.
    
    The action of $R$ on $T$ has non-trival free arc stabilizers, though it restricts to a free action
    of $F$ on $T$.  To see this note the two ways of foliating the rectangle with sides $\alpha , \beta$.
    In one foliation the leaves have gradient $1$, and in the other $-1$.
    In the orbit space $R\backslash X$ there are three rectangular discs attached to the same
    $1$-skeleton.  The foliation on each is one of the two foliations above and so two of them have
    the same foliation.  This means that there will be pairs of  points in the boundary which are joined
    by leaves in different $2$-cells and so there will be loops in leaves in the orbit space.
    These will not lift to loops in $X$ and so there will be non-trivial point stabilizers.
    The action of $R$ on the $\R$-tree will not be stable.
    It will restrict to the obvious free action of $F$ corresponding to a foliation by lines of gradient $1$ on 
    $E$.

   The tree $T$ has four orbits of branch points and the number of directions at a transversal
   of branch points will be $n_a, n_b, n_c$ and $n_d$.   If these integers are odd primes then 
   the only morphisms in the category of $\R$-trees with domain $T$ will be isomorphisms as there
   can be no folding at such a branch point.
   The tree corresponding to the foliation with dense leaves on the $2$-sphere with $4$ cone points will be another such $\R$-tree and there is a surjective homomorphism from the group of this orbifold to
   $R$.  However there is no morphism of $\R$-trees which includes this homomorphism as it is not
   an isomorphism.

    \bigskip
    \medskip
\noindent
{\bf Parallelepiped groups}
    
We now  generalize the construction by showing that there are $n$-dimensional parallelepiped groups, rectangle groups being the groups obtained when $n=2$.

We first consider the case $n=3$.
Consider the cube as shown in Fig 3.

Let $S$ be the group with  presentation

$$S = \langle a, b, c, d, e, f, g, h   \  |  \  ab^{-1} = cd^{-1} 
=ef^{-1} = gh^{-1}, ac^{-1} = bd^{-1} = eg^{-1} = fh^{-1}, $$
$$ ae^{-1} = bf^{-1} = cg^{-1} = dh^{-1} \rangle .$$

The form of this presentation suggests calling $S $ a {\it cube  } or even a {\it rectangular brick} group.
Note that there are $6$ rectangle subgroups corresponding to the $6$ faces of the cube $\langle a, b, c, d\rangle , \langle e, f, g, h \rangle ,$ $  \langle a, b, e, f \rangle ,$
$ \langle a, c, e, g \rangle , \langle b, d, f, h \rangle , \langle c, d, g, h \rangle$ and there are free abelian rank $3$ subgroups generated by 
$x = ab^{-1}, y= ac^{-1}, z = ae^{-1}$,  and $x' = b^{-1}a, y' = c^{-1}a, z' = e^{-1}a$.

The group $S$ has a decomposition as a free product with amalgamation

$$C =   \langle a, b, c, d  \rangle *_D \langle e, f, g, h  \rangle $$

where $D  = \langle  a{b}^{-1}, {b}^{-1}a, ac^{-1}, c^{-1}a  \rangle = \langle ef^{-1}, f^{-1}e, eg^{-1}, g^{-1}e \rangle $.
Here $ab^{-1} =x$ and $ac^{-1} = y$ commute and  $b^{-1}a = x'$ and $c^{-1}a = y'$ commute, and in fact, $D$ is a free product of two free abelian
rank $2$ subgroups.

There are two other such decompositions corresponding to the two vertical hyperplanes.

The vertex groups in the above decompositions are rectangle groups.

     \bigskip

 We construct
    a simply connected $3$-dimensional $S$-complex $X$ in which again there are $8$ orbits of vertices and edges
    and there is a sub-complex which is a three-dimensional Euclidean space  on which $x, y, z $ act by translation by $2$ units in the corresponding direction.  Let $J = \langle x, y, z\rangle$.
    The easiest way to construct  the space $X$ is as the cubing (see [Sa]) associated with the pairs
    $(C, D_1), (C,D_2), (C, D_3)$ where $D_1, D_2, D_3$ are the three splitting subgroups described
    above. The space $X$ will contain subcomplexes invariant under the rectangle groups associated with any two of these pairs and there will be three orbits of hyperplanes corresponding to
    the three splittings.

    Let $A, B, C, D, E, F, G, H$ be a tranversal for the vertex set of $X$, and let the stabilizers of the vertices be generated by $a, b, c, d, e, f, g, h$ respectively.

     Let $E_3$ be the subcomplex
    consisting of the $8$ $J$-orbits of a $2\times 2$-cube.  There are $20$
    $J$-orbits of $2$-cells and $8$ orbits of $3$-cells.  Under the action of $S$ two of the $8$ $3$-cells
    are in the same orbit and there are no other cases in which $J$-orbits are identified giving $7$
    $S$-orbits of $3$-cells in all.  The two $3$-cells that are in the same orbit only share a single vertex
    $A$ and are a $3$-cell and its translate by $a$.  Also under the action of $S$ six pairs of $2$-cells
    are identified giving $14$ $S$-orbits of $2$-cells.  If all the cubes containing a transversal of faces  lay in different
    orbits there would be $42$ orbits of $2$-cells.  In fact there are $14$.  Thus every face lies in $3$
    cubes.
    
    The above description - in which $S\backslash X$ is obtained from $J\backslash E_3$ by identifying one 
    pair
    of $3$-cells - also enables us to understand the decompositions corresponding to the hyperplanes
    geometrically.   In $J\backslash E_3$ a hyperplane gives two $2$-dimensional tori corresponding
    to $x = -1/2$ and $x = 1/2$.   Each of the tori is a $(2\times 2)$-rectangle in which opposite sides 
    are identified.   In $S\backslash X$ a $(1\times 1)$ rectangle in one torus is identified with 
    a $(1 \times 1)$- rectangle in the other torus.   The splitting group is the fundamental group
    of two tori with a rectangular disc identified, i.e. it is the free product of two free abelian groups 
    of rank two.
    
    As for rectangle groups we can determine all the simplicial $\R$-trees on which $S$ acts.
    Let $(p, q, r)$ be a triple of integers for which $1$ is the greatest common divisor.
    There is a surjective homomorphism $\theta : J \rightarrow {\bf Z}$ in which $\theta (x) = p, \theta (y) = 
    q$ and $\theta (z) = r$.  We assign lengths $p, q, r$ to the edges of $X$ in the orbits of $AB, AC, AE$
    respectively and foliate $X$ so that $E_3$ is foliated by planes given by $x + y + z = c$ for a constant 
    $c$. If we let $c$ range over all constants with fractional part $1/2$ then we obtain representatives of all non-parallel tracks in $E_3$.   In $J\backslash E_3$ these planes give $2$ non-parallel tracks, since $2$ is the greatest common divisor 
    of $2p, 2q, 2r$.   Each three cell of $E_3$ will intersect the transversal of non-parallel planes in $E_3$ in  
    $p+q+r -1$ distinct discs.  There is one such $3$-cell which intersects the planes for which
    $c$ takes the values $1/2, 3/2, \dots p+q+r - 1/2$.   In $S\backslash X$ the two tracks in $J\backslash
    E_3$ become a single track by identifying $p + q +r -1$ distinct pairs of discs.   The splitting group
    will be a free product of two free abelian rank two subgroups and a free group of rank
    $p +q +r -2$.

    We have discussed the case $n=2$ and $n=3$.   In general, corresponding to an $n$-dimensional
    parallelepiped there will be a group $P_n$ that has $n$ decompositions over a subgroup
    that is a free product of two rank $n-1$ free abelian groups.   Corresponding to these
    decompositions there is a $1$-connected $n$-dmensional cubing $X_n$.
    Let $E_n$ be  $n$-dimensional Euclidean space and let $J_n$ be  generated by translations of $2$ 
    units
    in each of the coordinate directions.  The space $P_n\backslash X_n$ is obtained from $J_n
    \backslash E_n$ by identifying an antipodal pair of $n$-cells in  $J_n\backslash E_n$ .
    
    Suppose now that we give give positive real lengths $\alpha _1, \alpha _2, \dots , \alpha _n$ of 
    translation to the generators 
    of $J_n$.  We choose these lengths so that they are independent over the rationals. We also foliate
     $E_n$ by $n-1$-dimensional affine subspaces which intersect
    each $2$-dimensional coordinate hyperplane in a straight line of gradient one.
    
    The leaves of this foliation are the points of an $\R$-tree on which $P_n$ acts with free arc 
    stabilizers.
    \bigskip
    \noindent

    \bigskip

    {\bf References}

\medskip
\noindent
\item{[BF]} M.Bestvina and M.Feighn, {\it Stable
actions of groups on real trees,} Invent. Math.
{\bf 121} (1995), 287 -321.
medskip
\noindent
\item{[BH]} M.R.Bridson and A. Haefliger, {\it Metric spaces of non-positive curvature}, \ \ Springer (1999).
\medskip
\noindent
\item{[DD] }W.Dicks and
M.J.Dunwoody, {\it Groups acting on graphs},
Cambridge University Press, 1989.

\medskip
\noindent
\item{[D]} M.J.Dunwoody, {\it Groups acting on real trees}, Extended Version (2006).

 http://www.maths.soton.ac.uk/pure/preprints.phtml.
 
 \medskip
\noindent
\item{[H]} A.Haefliger, {\it Complexes of groups and
orbihedra}, in Group theory from a geometric
viewpoint, Eds. Ghys,Haefliger, Verjovsky,
World Scientific 1991.
\medskip
\noindent
\item{[Sa]} M.Sageev, {\it Ends of group pairs and non-positively curved cube complexes},  Proc.
London Math. Soc. (3) {\bf  71} (1995) 585-617.

\medskip
\noindent
\item{[Sh]} P.B.Shalen,{\it Dendrology of groups: an introduction.}
Essays in Group Theory (edited by S.M.Gersten), MSRI Publication {\bf 8} (1987) 265-320.
\bigskip
Department of Mathematics

University of Southampton

Southampton

SO17 1BJ

M.J.Dunwoody@maths.soton.ac.uk
  \vfill
  \eject  
\input pictex

\input pictex
\midinsert
\line \bgroup \hss
\beginpicture
\setcoordinatesystem units <.071in, .1in>
\put {$ \ $} at 25 55
\setlinear 
\plot 0 0   0  50    /
\plot  10 0   10  50  /
\plot   20 0   20  50 /
\plot 30 0 30  50   /

\plot 40 0 40 50    /

\plot 50 0 50 50  /
\plot  -5 5   55 5  /
\plot  -5 15  55 15 /
\plot -5 25  55 25 /
\plot -5 35  55 35 /
\plot -5 45 55 45 /
\setshadegrid span <1pt>
\vshade 10 5 15 20 5 15 /
\vshade 30 5 15 40 5 15 /
\vshade 10 25 35  20 25 35 /
\vshade 0 15 25 10 15 25 /
\vshade 30 25 35 40 25 35 /

\vshade 20 15 25 30 15 25 /
\vshade 40 15 25 50 15 25 /
\vshade 0 35 45 10 35 45 /
\vshade 20 35 45 30 35 45 /
\vshade 40 35 45 50 35 45 /
\vshade 0 0 5  10 0  5 /
\vshade 20 0 5  30 0 5 /
\vshade 40 0 5 50 0 5  /
\vshade -5  5  15  0  5  15 /
\vshade -5  25  35  0 25 35 /
\vshade -5 45 50 0 45 50   /
\vshade  10 45  50 20 45 50   /
\vshade 30 45 50 40 45 50   /
\vshade 50  45 50 55 45 50   /
\vshade 50  25 35 55 25 35   /
\vshade 50  5 15 55 5 15   /
\put {$A$} at 19 24
\put {$B$} at  29 24
\put {$C$} at  19 16
\put {$D$} at  29  16
\put {$aB$} [r] at 9.5 24
\put {$cD$} [r] at 9.5 16
\put {$aC$} [r] at 19 34
\put {$bD$} [r] at 29 34
\put {$aD$} [r] at 9 34

\put {Fig 1} at  25 -5
\endpicture
\hss \egroup
\endinsert

\midinsert
\line \bgroup \hss
\beginpicture
\setcoordinatesystem units <.03in, .03in>
\put {$ \ $} at 25 55
\setlinear 
\plot 0 0   0  70    /
\plot  10 0   10  70  /
\plot   20 0   20  70 /
\plot 30 0 30  70   /

\plot 40 0 40 70    /

\plot 50 0 50 70  /
\plot  -5 5   55 5  /
\plot  -5 20  55 20 /
\plot -5 35  55 35 /
\plot -5 50  55 50 /
\plot -5 65 55 65 /

\plot -5 0   55 60 /
\plot  -5 2 55 62 /
\plot -5 4  55 64 /

\plot -5 6   55 66 /
\plot  -5 8 55 68 /

\plot  -5 10 55 70 /
\plot  -5 12 53 70 /
\plot -5 14  51 70 /

\plot -5 16   49 70 /
\plot  -5 18  47 70 /

\plot -5 20  45 70 /
\plot  -5 22 43 70 /
\plot -5 24   41 70 /

\plot -5 26 39 70 /
\plot -5  28 37 70 /
\plot -5 30  35 70  /
\plot -5 32 33 70 /
\plot -5 34 31 70 /
\plot -5 36 29 70 /
\plot  -5 38 27 70 /
\plot  -5 40 25 70 /
\plot  -5 42 23 70 /
\plot  -5 44 21 70 /
\plot  -5 46 19 70 /
\plot -5 48 17 70 /
\plot -5 50 15 70 /
\plot  -5 52 13  70 /
\plot  -5 54 11  70 /
\plot  -5 56 9 70 /
\plot  -5 58 7 70 /
\plot  -5 60 5 70 /
\plot  -5 62 3  70 /
\plot  -5 64 1  70 /
\plot  -5 66 -1 70 /
\plot  -5 68 -3 70 /
\plot -3 0 55 58 /
\plot -1 0 55 56 /
\plot  1 0 55 54 /

\plot 3 0 55 52 /
\plot 5 0 55 50 /
\plot 7 0 55 48 /
\plot 9 0 55 46 /
\plot 11 0 55 44 /
\plot 13 0 55 42 /
\plot 15 0 55 40 /
\plot 17 0 55 38 /
\plot 19 0 55 36 /
\plot 21 0 55 34 /
\plot 23 0 55 32 /
\plot 25 0 55 30 /
\plot 27 0 55 28 /
\plot 29 0 55 26 /
\plot 31 0 55 24 /
\plot 33 0 55 22 /
\plot 35 0 55 20 /
\plot 37 0 55 18 /
\plot 39 0 55 16 /
\plot 41 0 55 14 /
\plot 43 0 55 12 /
\plot 45 0 55 10 /
\plot 47 0 55 8 /
\plot 49 0 55 6 /
\plot 51 0 55 4 /
\plot 53 0 55 2 /
\plot 55 0 55 0 /

\put {Fig 2} at  25 -5
\endpicture
\hss \egroup
\endinsert
\bigskip
\vfill

\midinsert
\line \bgroup \hss
\beginpicture
\setcoordinatesystem units <.06in, .06in>
\put {$ \ $} at 25 55
\setlinear 
\plot 0 10   20 20  50  15 30 5 0 10    /
\plot  0 40   20  50  50 45 30 35 0 40  /
\plot   0 10   0 40  /
\plot  20 20 20 50    /

\plot  30 5 30 35     /

\plot 50 15 50 45   /

\put {$a$} at  -1 40
\put {$b$} at  20 51
\put {$d$} at  51 45
\put {$c$} at  29 34
\put {$e$}  at -1 10
\put {$f$}  at 20 19
\put {$h$}  at 51 15
\put {$g$}  at 30 4

\put {Fig 3} at  25 -5
\endpicture
\hss \egroup
\endinsert

\bye

\plot 0 0   0  50    /
\plot  10 0   10  50  /
\plot   20 0   20  50 /
\plot 30 0 30  50   /

\plot 40 0 40 50    /

\plot -5 0   55 42 /
\plot  -5 10 55 52 /
\plot -5 20 40 52 /
\plot -5 30 26.6  52 /
\plot -5 40 12.2 52 /
\plot 9.4 0 55 32 /
\plot 23.6 0 55 22 /
\plot 38 0 55 12 /

\bye

We want to show that there is a sequence of such spaces which tends to a limit
    where the limit corresponds to a psuedometric on $X$ in which points zero distance apart 
    are on a leaf of a foliation as shown in Fig 2.  In this diagram the plane is tessellated by
    rectangles with sides $\alpha,  \beta$ and the foliation is by lines with gradient one.
    Now $P$ is also generated by the stabilizers of the points $B, bA, C, cD$ which from
    Fig 1 can be seen to be the vertices of a parallelogram whose translates tessellate the $(x, y)$-plane.
    Let $X'$ be the graph with the same vertex set as $X$ but in which the edge set are the orbits 
    of the four sides of this parallelogram.
    We want to assign lengths to the edges so that it accords with the transverse distance on our
    foliation.   Clearly the length of the edges $(B, bA)$ and $(C, cD)$ will still be $\alpha $.
    Traversing the edge $(B, cD)$ corresponds to going $2\alpha $ in the direction $(B, aB)$ and $\beta $ in the reverse direction $(aB, cD)$.  Thus both $(B, cD)$ and $(bA, C)$  are assigned the length 
    $|\beta - 2\alpha |$.   There is a map $\theta : X(A, B, C, D, \alpha , \beta ) \rightarrow X(bA, B, C, cD,
    \alpha , |\beta - 2\alpha |)$ which is the identity map on the vertex set.   I thought originally that this would be a folding map if $\beta > \alpha $,
    but this is not necessarily the case.  It is possible that distances increase.
    Note that if $ \beta > \alpha $ then $|\beta - 2\alpha | < \beta$ and so the process of successively
    replacing the larger one of a pair of positive real numbers by the modulus of that number
    minus twice the smaller number will either terminate when one of the numbers is zero or
    both numbers will tend to zero.
    If one carries out the corresponding operations on the graphs $X(A, B, C, D, \alpha , \beta)$ then
    the distance between a pair of vertices  tends to a limit, which in the case of points in the $(x, y)$
    plane will be transverse measure on the foliation, i.e. $d((x_1, y_1), (x_2, y_2)) = |y_1-x_1 - y_2 +x_1|$.
This seems obvious but I am still trying to write down a proof.    
    

        If $\alpha , \beta $ are dependent over the rationals, then repeating this  operation
    will eventually terminate when one of the coefficients becomes zero.  This corresponds to
    one of the decompositions of $P$ as a free product with amalgamation discussed earlier.
    If $\alpha , \beta $ are independent then the limiting metric is that of an $\R$-tree.

 \bigskip
\noindent

Let $T$ be a simplicial $P_1$-tree.  We can assume that the distance between adjacent vertices
is $1$.
 Let $c = y_2y_1^{-1}, d =  y_3y_1^{-1} $.   These elements are of infinite order and 
commute.   If they fix a point of $T$, then so will $P_1$.  The other possibility
is that $c$, $d$ share an axis $A$.   Each of the points $y_1, y_2, y_3, y_4$ will 
stabilize a point of $A$. Suppose $y_i$ fixes the point $a_i, i =1,2,3,4$.   If $d(a_1, a_2) = \gamma $, then 
 $c$ has hyperbolic length $2\gamma $.  This means that $d(a_3, a_4) = \gamma, d(a_1, a_2)$
 and $d(a_2, a_4) = d(a_1, a_3) = \delta $.   This means that - possibly after relabelling - the points
 $a_1, a_2, a_3, a_4$ are arranged along $A$ in that order with the interval lengths $\gamma , \delta - \gamma, \gamma $.   Here $\gamma $ may be zero.   Either 
 
 Consider the case when $\gamma = 0$.  The stabilizer $X$ of $a_1 = a_2$ contains  $ \langle y_1, y_2 \rangle$  
 and the stabilizer $Y$ of $a_3 = a_4$ contains  $ \langle y_3, y_4 \rangle$.  The stabilizer $C$ of the edge $(a_1, a_3)$ is fixed by
 $c$.   
 Clearly $P_1 = X*_CB$.
 If say $y_1 \in C$ then $B \geq \langle y_1, y_3, y_4 \rangle = P_1$ and the decomposition is trivial.
Consider the word $y_1y_3^{-1}y_4y_2^{-1}$.  This is a relator in $P_1$ and so $y_3^{-1}y_4 \in C$
as $y_1 \in A -  C, y_2^{-1} \in B -C$.  But this means 
 which fixes a vertex, must fix the  incident edge.  Similarly $y_4^{-1}y_3$ must fix this edge
 and so it will be fixed by $y_3^2$ and $y_4^2$.   If either $y_3, y_4$ has odd order then
 there will be a vertex fixed by three of the generators.  But these will generate $P_1$ and so the action is trivial.

 \bye

 Conjugating $y_1$ by $c$ and leaving $y_2, y_3, y_4$ fixed, will translate  $a_1$ by $2\gamma$.
 If $\gamma > \delta - \gamma$,then the new points are arranged in the order $a_2, a_3, a_1, a_4$
 with distances $\delta - \gamma, 2\gamma - \delta, \delta - \gamma$.  If $\gamma \leq \delta - \gamma $
 then the new points are arranged $a_2, a_1, a_3, a_4$ with distances $\gamma,  \delta - 2\gamma, \gamma $.   In the first case, if $y_1' = cy_1^{-1}c^{-1} = y_2y_1^{-1}y_2^{-1}$,  then $y_1"y_2^{-1} = c = y_4y_3^{-1}$ and $y_3y_2'^{-1} = dc^{-1} = 
 y_4y_1'^{-1}$

There is an action of $P = P(n_1, n_2, n_3, n_4)$ on $T_1$ 
as follows:-  If $a \in A_1, b\in A_2, a+b \in A_3$

$$ P \rightarrow Isom (T_1),    x_1 \mapsto \xi _1,
x_2 \mapsto  \theta _{-a}\xi _2 \theta _a,  x_3 \mapsto \theta _{-b}\xi _3 \theta _b,
x_3 \mapsto \theta _{-a-b}\xi _4 \theta _{a+b}$$.

Let $G_1$ be the subgroup of $G$ generated by $2A$ and the elements $x_1, x_2, x_3, x_4$.
We show that $T$ is a $G_1$-tree, by constructing isometries for each generator.
For $x_1$,  let  $x_1w$ be the word $x_1w$ if $a_1 \not= 0$ or if $a_1 = 0$ and $f_1 \not = x^{-1}$.
If $a_1 = 0$ and $ f_1 = x_1^{-1}$ then $x_1w = a_2f_2 \dots a_{n-1}f_{n-1}a_n$.
We now construct an  isometry corresponding to $z \in 2A$.  This is to be a hyperbolic isometry with axis $R$ where
$R$ is the union of the points $w = a_1$ for some $a_1 \in A, a_1 \geq 0$ and the points $x_1^{-1}a_2$ for $a_2 \in A, a_2 \geq 0$.   Thus $zw = a_1'f_1a_2 f_2 \dots a_{n-1}f_{n-1}a_n$ where $a_1' = a_1+z$ if $a_1 + z > 0$ and
$zw = x_1a_1'x_{i_1}^{-1}f_1a_2 f_2 \dots a_{n-1}f_{n-1}a_n$ where $a_1' = -a_1 - z$ if $a_1 + z < 0$.
Here the inclusion of $x_{i_1}^{-1}$ in the word reflects the fact that the relevant branch point is approached
in a different direction, since $a_1$ and $a_1+z$ have different signs.   If $f_1 = x_{i_1}$ then
we get a shorter word by cancellation in which $a_1' =-a_1 - z + a_2$.
If $a_1+ z = 0$ then $a_1 \in 2A = A_1$ and so $f_1 \in \langle x_1 \rangle$ and we put
$zw = f_1a_2 f_2 \dots a_{n-1}f_{n-1}a_n$.  
Thus if the natural direction to leave a point approached along $d$ is $d'$ then the natural direction
to leave it approaching along $d'$ is $x_id'$.

Consider the action of $g = zx_1^{-1}$ on $T$.
If $a_1 > z $ then $gw = (a_1- z) f_1 \dots a_{n-1}f_{n-1}a_n$, and if $0 \leq  a_1 \leq z$ 
then $gw  = (z - a_1)x_{i_1}^{-1} f_1 \dots a_{n-1}f_{n-1}a_n$.  Thus the point 
$q = z/2$ is fixed, and the directions at $d$ are rotated according to $x_{i_1}^{-1}$.
For $q = a, b, a + b$ respectively we get actions of $x_2, x_3 , x_4$ on $T$ via $(2q)x_1^{-1}$. 
Note that $x_2^{-1}x_1$ and $x_4^{-1}x_3$ both induce a translation of length $b$ along $\R$
and in fact the same action on $T$.
Thus the orbifold group
$$G = \langle x_1, x_2, x_3, x_4 | x_1^{n_1} = x_2^{n_2} = x_3^{n_3} = x_4^{n_4} =1.
x_2^{-1}x_1 = x_4^{-1}x_3 \rangle $$
has an action on $T$.
However this action is not faithful
since the translations $x_3^{-1}x_1$ and $x_4^{-1}x_2$ are also equal and commute with
$x_2^{-1}x_1$.

[D1] M.J.Dunwoody, {\it Groups acting on protrees,} J. London
Math. Soc. {\bf 56} (1997) 125-136.
\medskip
\noindent
[D2] M.J.Dunwoody, {\it Folding sequences,} Geometry and
Topology Monographs {\bf 1} (1998) 143-162.
\medskip
\noindent
[D3] M.J.Dunwoody, {\it A small unstable action on a tree,}
Math. Research Letters {\bf 6} (1999) 697-710. 
\medskip
\noindent
[D4] M.J.Dunwoody, {\it A tessellation and a
group acting on trees,}

http://www.maths.soton.ac.uk/pure/preprints.phtml.
\medskip
\noindent
[DD] W.Dicks and
M.J.Dunwoody, {\it Groups acting on graphs},
Cambridge University Press, 1989.
\medskip
\noindent
[H] A.Haefliger, {\it Complexes of groups and
orbihedra}, in Group theory from a geometric
viewpoint, Eds. Ghys,Haefliger, Verjovsky,
World Scientific 1991.
\medskip
\noindent
[K] P.H.Kropholler, {\it A note on centrality
in $3$-manifold groups,} Math. Proc. Camb.
Phil. Soc. {\bf 107} (1990) 261-266.
\medskip
\noindent
[L] G.Levitt, {\it Non-nesting actions on real trees,} Bull. London Math. Soc.
{\bf 30}, (1998) 46-54.
\medskip
\noindent
[LP] G.Levitt and F.Paulin, {\it
Geometric group actions on trees}, Amer. J.
Math {\bf 119} (1997) 83-102.
\medskip
\noindent
[RS] E.Rips and Z.Sela, {\it Cyclic splittings
of finitely presented groups and the canonical
JSJ-decomposition}, Annals of Maths. {\bf 146}
(1997) 53-104.
\medskip
\noindent
[S] R.Skora, {\it Splittings of surfaces,} J. Amer. Math.Soc. 9 (1996) 605-616.

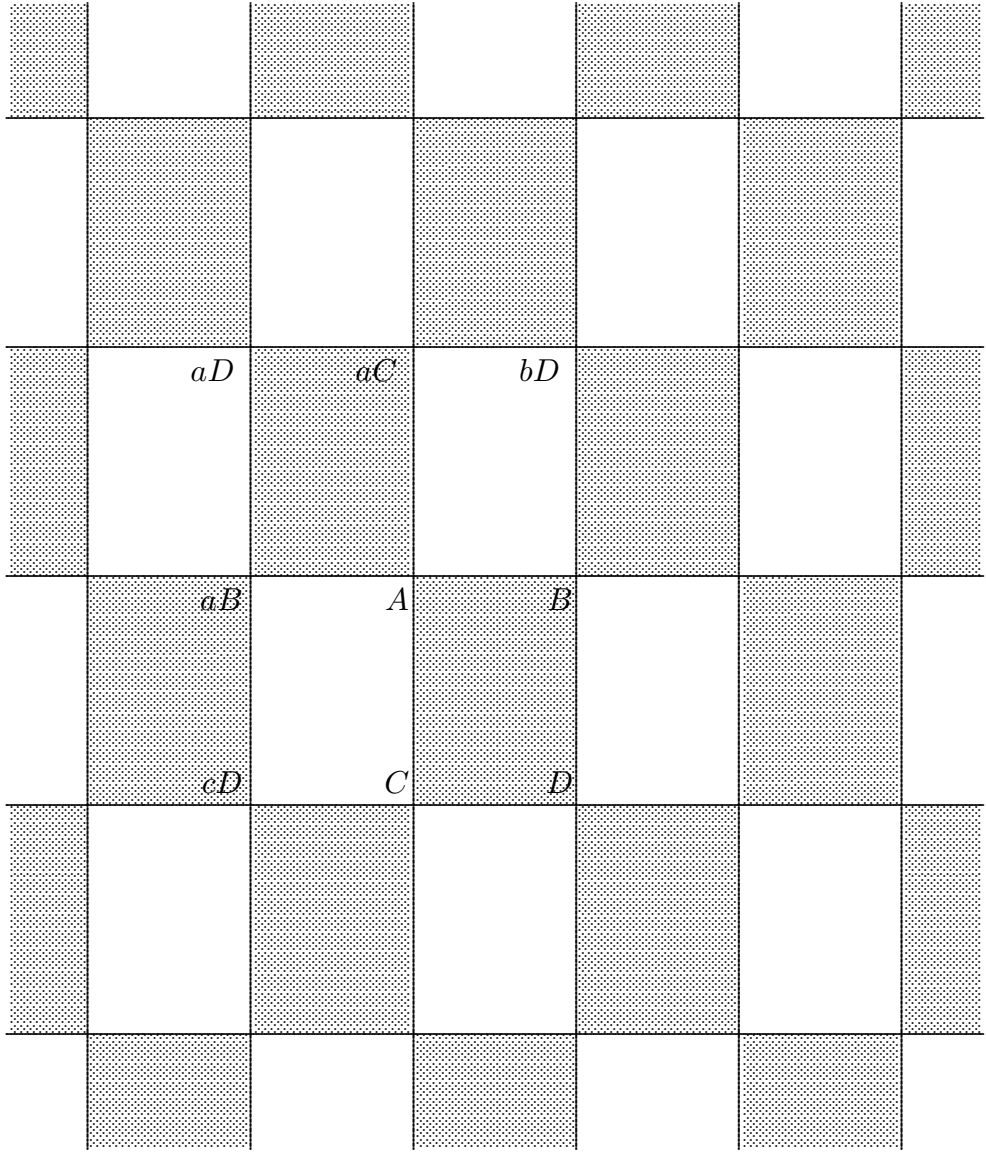
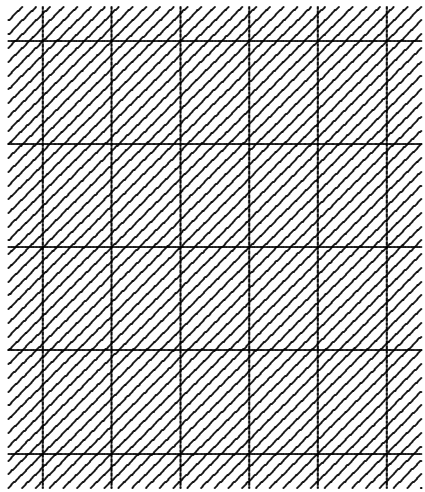
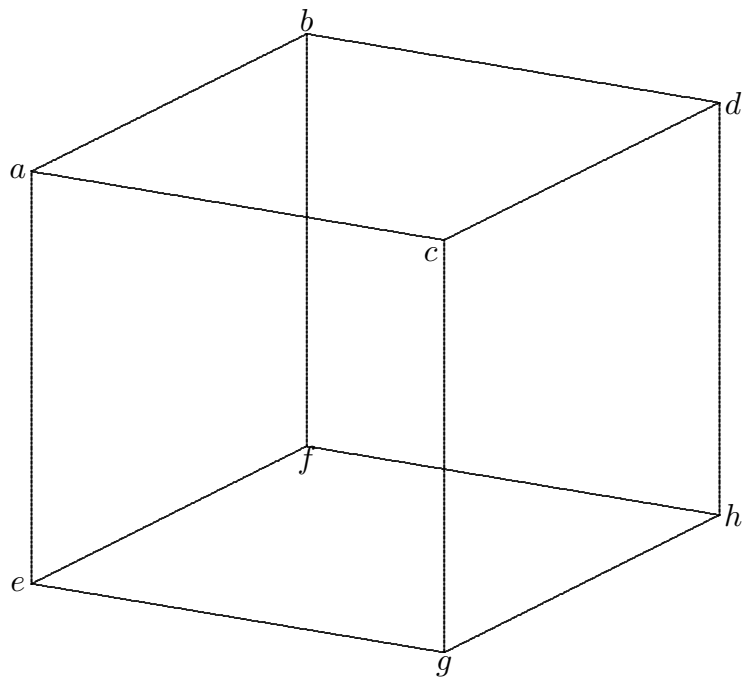
\bye